\documentclass[12pt]{amsart}
\usepackage{amsfonts,amssymb,amscd,amsmath,enumerate,verbatim,calc,latexsym,pstcol,pst-plot}
\def\frk{\mathfrak}               

\def\Phi{{\frk N}}
\def\opn#1#2{\def#1{\operatorname{#2}}} 
\opn\chara{char} 
\opn\length{\ell} 
\opn\pd{pd} 
\opn\rk{rk}
\opn\projdim{proj\,dim} 
\opn\injdim{inj\,dim} 
\opn\rank{rank}
\opn\depth{depth} 
\opn\grade{grade} 
\opn\height{height}
\opn\embdim{emb\,dim} 
\opn\codim{codim}

\opn\Tr{Tr} 
\opn\bigrank{big\,rank}
\opn\superheight{superheight}
\opn\lcm{lcm}
\opn\trdeg{tr\,deg}
\opn\reg{reg} 
\opn\lreg{lreg} 
\opn\ini{in} 
\opn\lpd{lpd}
\opn\size{size}
\opn\mult{mult}
\opn\dist{dist}
\opn\cone{cone}
\opn\lex{lex}
\opn\rev{rev}
\opn\div{div} \opn\Div{Div} \opn\cl{cl} \opn\Cl{Cl}
\opn\Spec{Spec} \opn\Supp{Supp} \opn\supp{supp} \opn\Sing{Sing}
\opn\Ass{Ass} \opn\Min{Min}
\opn\Ann{Ann} \opn\Rad{Rad} \opn\Soc{Soc}
\opn\Syz{Syz} \opn\Im{Im} \opn\Ker{Ker} \opn\Coker{Coker}
\opn\Am{Am} \opn\Hom{Hom} \opn\Tor{Tor} \opn\Ext{Ext}
\opn\End{End} \opn\Aut{Aut} \opn\id{id} \opn\ini{in}

\opn\nat{nat}
\opn\pff{pf}
\opn\Pf{Pf} \opn\GL{GL} \opn\SL{SL} \opn\mod{mod} \opn\ord{ord}
\opn\Gin{Gin}
\opn\Hilb{Hilb}\opn\adeg{adeg}\opn\std{std}\opn\ip{infpt}
\opn\Pol{Pol}
\opn\sat{sat}
\opn\Var{Var}
\opn\Gen{Gen}

\opn\aff{aff} \opn\con{conv} \opn\relint{relint} \opn\st{st}
\opn\lk{lk} \opn\cn{cn} \opn\core{core} \opn\vol{vol}
\opn\link{link} \opn\star{star}
\opn\gr{gr}

\def\pot#1#2{#1[\kern-0.28ex[#2]\kern-0.28ex]}

\opn\dirlim{\underrightarrow{\lim}}
\opn\inivlim{\underleftarrow{\lim}}
\def\Implies{\ifmmode\Longrightarrow \else
        \unskip${}\Longrightarrow{}$\ignorespaces\fi}
\def\implies{\ifmmode\Rightarrow \else
        \unskip${}\Rightarrow{}$\ignorespaces\fi}
\def\iff{\ifmmode\Longleftrightarrow \else
        \unskip${}\Longleftrightarrow{}$\ignorespaces\fi}

\let\:=\colon
\newtheorem{Theorem}{Theorem}[section]
\newtheorem{Lemma}[Theorem]{Lemma}
\newtheorem{Corollary}[Theorem]{Corollary}

\newtheorem{Remark}[Theorem]{Remark}

\newtheorem{Example}[Theorem]{Example}

\newtheorem{Conjecture}[Theorem]{Conjecture}

\let\epsilon\varepsilon
\let\phi=\varphi
\let\kappa=\varkappa
\def\qed{\ifhmode\textqed\fi
      \ifmmode\ifinner\quad\qedsymbol\else\dispqed\fi\fi}
\def\textqed{\unskip\nobreak\penalty50
       \hskip2em\hbox{}\nobreak\hfil\qedsymbol
       \parfillskip=0pt \finalhyphendemerits=0}
\def\dispqed{\rlap{\qquad\qedsymbol}}

\opn\dis{dis}
\def\pnt{{\raise0.5mm\hbox{\large\bf.}}}

\opn\Lex{Lex}

\begin{document}
\title{Betti numbers of chordal graphs and $f$-vectors
of simplicial complexes}
\author{Takayuki Hibi, Kyouko Kimura and Satoshi Murai}
\thanks{
{\bf 2000 Mathematics Subject Classification:}
Primary 13P99, 13F55; Secondary 52B05. \\
\, \, \, {\bf Keywords:}
monomial ideal, Betti sequence, simplicial complex, $f$-vector, 
chordal graph.
\\
\, \, \, 
The third author is supported by JSPS Research
Fellowships for Young Scientists. 
}
\address{Takayuki Hibi,
Department of Pure and Applied Mathematics,
Graduate School of Information Science and Technology,
Osaka University,
Toyonaka, Osaka 560-0043, Japan} 
\email{hibi@math.sci.osaka-u.ac.jp}
\address{Kyouko Kimura,
Department of Pure and Applied Mathematics,
Graduate School of Information Science and Technology,
Osaka University,
Toyonaka, Osaka 560-0043, Japan}
\email{kimura@math.sci.osaka-u.ac.jp}
\address{Satoshi Murai,
Department of Mathematics,
Graduate School of Science,
Kyoto University,
Sakyo-ku, Kyoto 606-8502, Japan}
\email{murai@math.kyoto-u.ac.jp}
\begin{abstract}
Let $G$ be a chordal graph and $I(G)$ its edge ideal. 
Let
$\beta(I(G)) = (\beta_0, \beta_1, \ldots, \beta_p)$
denote  
the Betti sequence of $I(G)$, where $\beta_i$ stands for 
the $i$th total Betti number of $I(G)$ and where $p$ is
the projective dimension of $I(G)$.  
It will be shown that there exists 
a simplicial complex $\Delta$ of dimension $p$
whose $f$-vector $f(\Delta) = (f_0, f_1, \ldots, f_{p})$ 
coincides with $\beta(I(G))$.
\end{abstract}
\maketitle

\section*{Introduction}
Let $S=K[x_1,\dots,x_n]$ be the polynomial ring 
in $n$ variables over a field $K$
with each $\deg x_i = 1$.
The {\em Betti sequence} of a homogeneous ideal 
$I \subset S$ is
the sequence
\[
\beta(I)=(\beta_0(I),\beta_1(I),\dots,\beta_p(I)),
\]
where each $\beta_i(I)$ stands for 
the $i$th total Betti number of $I$ and 
where $p = \projdim I$ is the projective dimension of $I$.
One has
$\sum_{i=-1}^{p} (-1)^{i} \beta_i(I) = 0$ with $\beta_{-1}(I) = 1$.

Let $\Delta$ be a simplicial complex and 
\[
f(\Delta) = (f_0, f_1, \ldots, f_{d-1})
\] its {\em $f$-vector},
where each $f_i = f_i(\Delta)$ stands for the number of faces
of $\Delta$ of dimension $i$ 
and where $d - 1$ is the dimension $\Delta$.
Recall that $\Delta$ is {\em acyclic} (over $K$)
if its reduced homology group
${\tilde H}_i(\Delta;K)$ 
with coefficients $K$ 
vanishes for all $i$. 
Thus in particular if $\Delta$ is acyclic,
then its $f$-vector 
satisfies $\sum_{i=-1}^{d-1} (-1)^i f_i = 0$
with $f_{-1} = 1$.
 
Peeva and Velasco \cite{PV} succeeded in proving that,
given an acyclic simplicial complex $\Delta$, 
there exists a monomial ideal $I$ whose Betti sequence
$\beta(I)$ coincides with the $f$-vector $f(\Delta)$.
In general, the converse is, however, false.
Let $n = 6$ and
$I = (x_1x_2, x_2x_3, x_3x_4, x_4x_5, x_5x_6, x_1x_6)$.
Then $\dim S / I = 3$, $\depth S / I = 2$
and $p = 4$.  One has $\beta(I) = (6, 9, 6, 2)$.
If a simplicial complex $\Delta$ 
possesses $2$ faces of dimension $3$, then 
$\Delta$ possesses at least $7$ faces 
of dimension $2$.
It then follows that 
there exists {\em no} simplicial complex
$\Delta$ of dimension $3$ 
with $(6, 9, 6, 2)$ its $f$-vector.

On the other hand, in Example \ref{nonacyclic},
one can find a Cohen--Macaulay monomial ideal
$I$, i.e., $S / I$ is a Cohen--Macaulay ring, 
whose Betti sequence is the $f$-vector of
a simplicial complex, but {\em not} the $f$-vector
of an acyclic simplicial complex. 

It is natural to ask which monomial ideals $I$
enjoy the property that there exists
a simplicial complex (or acyclic simplicial
complex) $\Delta$ whose $f$-vector coincides
with the Betti sequence of $I$.
The purpose of the present paper is to establish
the research project
on finding a natural class $\mathcal C$ 
of monomial ideals such that,
for each ideal $I$ belonging to $\mathcal C$, 
the Betti sequence $\beta(I)$ is the $f$-vector
of a simplicial (or an acyclic simplicial) complex.



First, in Section $1$, we summarize
several answers, which are easily or directly 
obtained from well-known facts.
The topics discussed will include
monomial ideals with small projective dimensions,
cellular resolutions, 
componentwise linear ideals and
pure resolutions.

Now, Section $2$ is the highlight of this paper.
Let $G$ be a finite graph on the vertex set $V$
and $E(G)$ the edge set of $G$. 
We write $S = K[\{ x : x \in V \}]$ for the polynomial ring
in $|V|$ variables over a field $K$ 
with each $\deg x = 1$.  The edge ideal 
of $G$ is the ideal $I(G)$ of $S$ generated by
those monomials $x y$ with $\{x, y\} \in E(G)$.  
Recall that a finite graph $G$ is chordal
if each cycle of $G$ of length $> 3$ has a chord.
Theorem \ref{claim:chordal_graph}
guarantees that, for an arbitrary chordal graph $G$, 
there exists 
a simplicial complex $\Delta$ 
whose $f$-vector 
coincides with $\beta(I(G))$.
The recursive-type formula 
due to H\`{a} and Van Tuyl
\cite{Ha-Tuyl} 
will be indispensable to achieve the proof of 
Theorem \ref{claim:chordal_graph}.

Finally, in Section $3$, we study Gorenstein monomial ideals.
It follows that 
the Betti sequence of a Gorenstein monomial ideal
$I$ with $\projdim (I) \leq 3$ is the $f$-vector
of an acyclic simplicial complex.
On the other hand, we can characterize the possible Betti numbers 
of Gorenstein monomial ideals $I$ with $\projdim (I) = 3$. 
Moreover, it will be proved that, given
integers $m \geq 4$ and $p \geq 3$,
there exists a Gorenstein monomial ideal $I$ of
$K[x_1, \ldots, x_n]$, where $n$ is enough large,
with $\beta_0(I) = m$ 
and $\projdim(I) = p$ if and only if  
$m \geq p + 1$ with $m \ne p+2$.

\section{Betti sequences and 
acyclic simplicial complexes}

The present section is
a summary of several answers, which are easily or directly 
obtained from well-known facts, for 
the problem of finding a natural class $\mathcal C$ 
of monomial ideals such that,
for each ideal $I$ belonging to $\mathcal C$, 
the Betti sequence $\beta(I)$ is the $f$-vector
of a simplicial (or an acyclic simplicial) complex.

First, recall a combinatorial characterization 
of $f$-vectors of acyclic simplicial complexes 
due to Gil Kalai \cite{Ka}.

\begin{Lemma}[Kalai]
\label{Kalaiacyclic}
A vector $f=(f_0,f_1,\dots,f_{d-1})$ of positive integers
is the $f$-vector of
an acyclic simplicial complex
of dimensional $d - 1$ 
if and only if there exists
a simplicial complex $\Delta'$
of dimension $d - 2$
with $f(\Delta') = (f'_0, f'_1, \ldots, f'_{d-2})$
such that $f_i=f'_i + f'_{i-1}$ for all $i$,
where $f'_{-1}=1$ and
$f'_{d-1} = 0$.
\end{Lemma}


\noindent
{\bf (1.1) Monomial ideals with small projective dimensions}

Let $I \subset S$ be a monomial ideal with 
$\projdim(I)\leq 2$
and $\beta(I)=(n,\beta_1,\beta_2)$.
One has $1 - n + \beta_1 -\beta_2=0$.
It follows from the Taylor resolution of $I$ that
there exists an integer $c \geq 0$
such that $\beta_1= {n \choose 2} -c$.
Thus $\beta(I)=(n,{n \choose 2} -c, {n-1 \choose 2} -c)$.
Since $(n-1,{n-1 \choose 2} - c)$ is the $f$-vector of a simplicial complex,
Lemma \ref{Kalaiacyclic} says that
$\beta(I)=(n, {n \choose 2} -c, {n-1 \choose 2} -c)$
is the $f$-vector of an acyclic simplicial complex.

\begin{Theorem}
Let $I \subset S$ be a monomial ideal
with $\projdim(I)\leq 2$.
Then $\beta(I)$
is the $f$-vector of an acyclic simplicial complex.
\end{Theorem}

%


\noindent
{\bf (1.2) Cellular resolutions}

The cellular resolution was 
introduced by Bayer and Sturmfels \cite{BS}.
Let $I \subset S$ be a monomial ideal
and $\mathbf F_\bullet$ a $\mathbb{Z}^n$-graded
free resolution of $S/I$.
The complex $\mathbf F_\bullet \otimes_S S/(x_1-1,\dots,x_n-1)$
of $K$-vector spaces is called the \textit{frame} of $\mathbf{F}_\bullet$.
We say that $\mathbf F_\bullet$ is supported by a CW-complex $\Delta$
if its frame is equal to the augmented oriented chain complex of $\Delta$.
If a free resolution is supported by a CW-complex $\Delta$, then
$\Delta$ must be acyclic (\cite[Proposition 1.2]{BS}).
Thus if a minimal free resolution is supported by a simplicial complex, 
then its Betti sequence must be the $f$-vector 
of an acyclic simplicial complex.

A monomial ideal $I \subset S$ is said to be \textit{generic} if,
for all pairs of generators $u=x_1^{a_1}\dots x_n^{a_n}$
and $v=x_1^{b_1}\dots x_n^{b_n}$ of $I$,
one has $a_k\ne b_k$ or $a_k=b_k=0$ for all $k$.
It was proved by Bayer, Peeva and Sturmfels \cite{BPS} that
a generic monomial ideal has  
a minimal free resolution which is supported by a simplicial complex.

\begin{Theorem}
Let $I$ be a generic monomial ideal.
Then $\beta(I)$ is the $f$-vector of an acyclic simplicial complex.
\end{Theorem}

We say that a CW-complex $\Delta$ satisfies the \textit{intersection property}
if the intersection of two faces of $\Delta$ is again a face of $\Delta$.
For example, 
all simplicial complexes as well as all polyhedral complexes
satisfy the intersection property.
Bj\"oner and Kalai \cite{BK2} proved that
if $\Delta$ is an acyclic CW-complex satisfying the intersection property,
then the $f$-vector of $\Delta$ is the $f$-vector of an acyclic simplicial complex.

\begin{Theorem}
Suppose that 
the minimal free resolution of a monomial ideal $I \subset S$ 
is supported by a CW-complex
satisfying the intersection property. 
Then $\beta(I)$ is the $f$-vector of an acyclic simplicial complex.
\end{Theorem}

Velasco \cite{Ve} studied 
minimal free resolutions 
which are not supported by a CW-complex
by means of
the nearly scarf ideal introduced in \cite{PV}.
Let $\Omega$ be a simplicial complex with the vertex set 
$[n]=\{1,2,\dots,n\}$
which is not the boundary of a simplex.
The \textit{nearly scarf ideal} $J_\Omega$ of $\Omega$ is the monomial ideal
of the polynomial ring $K[x_\sigma: \sigma \in \Omega \setminus \{\emptyset\}]$
generated by 
$
\{
\prod_{\sigma \in \Omega,\ v \not \in \sigma } x_\sigma : v \in [n]
\}.$
It is known \cite{PV} that the graded Betti numbers of $J_\Omega$
is given by
$$\beta_i(J_\Omega) = f_i(\Omega) + \dim_K \tilde H_{i-1}(\Omega;K),
\, \, \, \, \, \, \, \, \, \, i \geq 0.$$
On the other hand, 
Bj\"orner--Kalai Theorem (\cite{BK1}),
which gives a characterization of the $(f,\beta)$-pairs 
of simplicial complexes,
guarantees that,
for an arbitrary simplicial complex $\Delta$
with $f(\Delta) = (f_0, f_1, \ldots, f_{d-2})$,
the vector $(f_0',\dots,f'_{d-1})$
defined by setting 
$f_i'=f_i + \dim_K \tilde H_{i-1}(\Delta;K)$
is the $f$-vector of an acyclic simplicial complex.

\begin{Theorem}
Let $J_\Omega$ be the nearly scarf ideal of $\Omega$.
Then $\beta(J_\Omega)$ is the $f$-vector of an acyclic simplicial complex.
\end{Theorem}


\noindent
{\bf (1.3) Componentwise linear ideals}

One of the most famous classes of monomial ideals for which
the formula of graded Betti numbers is known is the class of stable ideals.
Recall that a monomial ideal $I \subset S$ is \textit{stable}
if, for all monomials $u  \in I$ and for all $1 \leq i < m(u)$, 
one has $u x_i/ x_{m(u)} \in I$,
where $m(u)$ is the maximal integer $k$ such that $x_k$ divides $u$.
Let $I$ be a stable ideal and $G(I)$ the minimal set of 
monomial generators of $I$.
Write $m_k(I)$ for the number of monomials 
$u \in G(I)$ with $m(u)=k$.
Eliahou and Kervaire \cite{EK} proved that
\[
\beta_i(I)
= \sum_{k=i+1}^n m_k(I) {k-1 \choose i} 
\]
for all $i \geq 0$.

A homogeneous ideal $I \subset S$ is said to have a \textit{$k$-linear resolution} if $\beta_{i,i+j}(I)=0$ whenever $j \ne k$.
A homogeneous ideal $I \subset S$ is said to be 
\textit{componentwise linear} (\cite{HH})
if, for all integers $k \geq 0$, the ideal $I_{\langle k \rangle}$
which is generated by the homogeneous polynomials of degree $k$
belonging to $I$ 
has a $k$-linear resolution.
A \textit{quasi-forest} is a simplicial complex $\Delta$
whose Stanley--Reisner ideal $I_\Delta$ has a $2$-linear resolution.
It is known (Fr\"oberg \cite{F}) that 
a quasi-forest is the clique complex of
a chordal graph.

\begin{Theorem}
Let $\beta=(\beta_0,\beta_1,\dots,\beta_p)$ with $p \leq n-1$ 
be a sequence of integers.
The following conditions are equivalent:
\begin{itemize}
\item[(i)] There exists a componentwise linear ideal $I \subset S$ 
with $\projdim(I) =p$ such that $\beta(I)=\beta$;
\item[(ii)] There exists a stable ideal $I \subset S$ with 
$\projdim(I) =p$
such that $\beta(I)=\beta$;
\item[(iii)]
There exists a sequence $c_1,\dots,c_{p+1}$ of positive integers
with $c_1=1$ such that 
$\beta_i= \sum_{k=1}^{p+1}c_k {k-1 \choose i}$ for all $i \geq 0$;
\item[(iv)] There exists an acyclic quasi-forest $\Delta$ of dimension $p$ 
such that $\beta=f(\Delta)$.
\end{itemize}
\end{Theorem}

\begin{proof}
First, (i) $\Leftrightarrow$ (ii) is known (\cite[Lemma 1.4]{CHH}).
Second,
(ii) $\Rightarrow$ (iii) follows from Eliahou--Kervaire formula
and the fact that if $I$ is a stable ideal and $m_k(I) \ne 0$
for some $k>0$, then $m_\ell(I)\ne 0$ for all $1 \leq \ell <k$ 
(\cite[Lemma 1.3]{HHMT}).
Third, to prove (iii) $\Rightarrow$ (ii), 
we introduce the monomial ideal $I$
generated by
$$\bigcup_{i=1}^{p+1} \big\{(x_1^{c_2} \cdots x_{i-2}^{c_{i-1}})x_{i-1}^{c_i+1 -k} x_i^{c_{i+1}+k} : k=1,\dots,c_i\big\},$$
where $c_{p+2}=0$.
It follows that $I$ is stable and 
$(m_1(I), \ldots, m_{p+1}(I))=(c_1,\dots,c_{p+1})$.

Finally, (iii) $\Leftrightarrow$ (iv) will be shown.
It is known \cite{HHMTZ} that
$f=(f_0,f_1,\dots,f_{p-1})$ is the $f$-vector of a quasi-forest 
of dimension $p - 1$ if and only if
there exists a sequence of positive integers $b_1,\dots,b_p$
such that $f_{i-1}= \sum_{k=1}^p b_k {k-1 \choose i-1}$ for all $i \geq 1$.
If $\Delta$ is a quasi-forest, then 
it follows from \cite[Theorem 7.1]{He} that 
its algebraic shifted complex $\Sigma$ is again a quasi-forest.
If $\Delta$ is acyclic then $\Sigma$ must be a cone (\cite{Ka}).
However, if a quasi-forest $\Sigma$ is a cone, then it must be a cone
of a quasi-forest.
These facts guarantee that
$f=(f_0,f_1,\dots,f_p)$ is the $f$-vector of an acyclic 
quasi-forest of dimension $p$ if and only if
$f$ is the $f$-vector of a cone of a quasi-forest
of dimension $p - 1$.
The latter condition is equivalent to saying that
there exists a sequence of positive integers $b_1,\dots,b_p$
such that $f_{i-1}= \sum_{k=1}^{p} \{b_k {k-1 \choose i-1} + b_k{k-1 \choose i-2}\}=\sum_{k=1}^{p} {b_k} {k \choose i-1}$ for all $i \geq 2$ and $f_{0}=1+\sum_{k=1}^p b_k$.
Set $c_1=1$ and $c_k=b_{k-1}$ for $k=2,3,\dots,p+1$.
Then the sequence $c_1,\dots,c_{p+1}$ satisfies the conditions of (iii),
as desired.
\, \, \, \, 
\end{proof}


\noindent
{\bf (1.4) Pure resolutions}

We discuss the question whether Betti sequences of monomial ideals 
with pure resolutions are 
$f$-vectors of simplicial complexes. 
We say that a homogeneous ideal $I \subset S$ 
has a \textit{pure resolution} if 
its minimal free resolution is of the form 
\begin{displaymath}
  0 \longrightarrow  S(-c_p)^{\beta_p}
    \longrightarrow  S(-c_{p-1})^{\beta_{p-1}}
    \longrightarrow \cdots 
    \longrightarrow  S(-c_{0})^{\beta_{0}}
    \longrightarrow I 
    \longrightarrow 0. 
\end{displaymath}

Let $v > d \geq 1$ and
$C(v,d)$ the cyclic polytope \cite[p.\ 59]{Stanley}
of dimension $d$ with $v$ vertices.
Since $C(v,d)$ a simplicial polytope,
its boundary $\partial C(v,d)$ defines
a simplicial complex $\Delta(C(v,d))$,
called the boundary complex
of $C(v,d)$.  
It is known \cite[Proposition 3.1]{TH}
that, when $d$ is even, the Stanley--Reisner ideal
$I_{\Delta(C(v,d))}$ 
(\cite[p.\ 53]{Stanley})
of $\Delta(C(v,d))$
has a pure resolution.

\begin{Example}
{\em
Let $v = 7$ and $d = 2$.
Then the Betti sequence of $I_{\Delta(C(7,2))}$ is 
$(14,35,35,14,1)$.
In particular $(14,35,35,14,1)$ 
is the Betti sequence 
arising from a pure resolution.  However, 
it turns out that
$(14,35,35,14,1)$ 
cannot be the Betti sequence 
arising from a linear resolution.
}
\end{Example}

\begin{Example}
\label{nonacyclic}
{\em
In \cite{BrHi95} it is shown that there exists a simplicial
complex $\Delta$ such that 
(i)
$I_\Delta$ has a pure, but not a linear resolution;
(ii)
the Betti sequence of $I_{\Delta}$ is
$\beta(I_{\Delta}) = (14,21,14,6)$;
(iii)
the Stanley--Reisner ring 
$K[\Delta] = S / I_\Delta$ 
(\cite[p. 53]{Stanley}) 
is Cohen--Macaulay.  
Now, Kruskal--Katona theorem 
\cite[p. 55]{Stanley} says that 
$(14,21,14,6)$ is the $f$-vector
of a simplicial complex.  However,
by using Lemma \ref{Kalaiacyclic}
it turns out that 
$(14,21,14,6)$ cannot be the $f$-vector
of an acyclic simplicial complex.  
}
\end{Example}

\begin{Theorem}
  \label{claim:cycle}
  If $d$ is even, then the Betti sequence of $I_{\Delta (C(v,d))}$ is 
  the $f$-vector of a simplicial complex. 
\end{Theorem}

\begin{proof}
  Let $d=2d'$
  and 
  $\beta(I_{\Delta (C(v,d))}) = (\beta_0, \ldots, \beta_{v-2d'-1})$.
  It follows from \cite{TH} that
  \begin{eqnarray}
  \label{cyclic}
  \beta_i
  = \binom{v-d'-1}{d'+i+1} \binom{d'+i}{d'}
       + \binom{v-d'-1}{i} \binom{v-d'-i-2}{d'}
  \end{eqnarray}
  for $i < v-2d'-1$ and $\beta_{v-2d'-1} = 1$.

  \par
  Let $v=d+1$.  Then the Betti sequence of 
  $I_{\Delta (C(v,d))}$ is $(1)$, which is the $f$-vector 
  of a $0$-simplex. Let $v \geq d+2$.
Our proof will be done by using induction on $d'$. 

Let $d'=1$. Then
  $\Delta (C(v,2))$ is a cycle with $v$ vertices. 
  We show that, by using induction on $v$, the Betti 
  sequence $\beta(I_{\Delta(C(v,2))})$ is the $f$-vector 
  of a simplicial complex. 
  When $v=4$, the Betti sequence of $I_{\Delta (C(v,2))}$
  is $(2,1)$, which is the $f$-vector of a $1$-simplex. 
  Let $v > 4$ and suppose that
  there exists a simplicial complex $\Gamma (v-1)$ 
  such that $f (\Gamma (v-1)) = \beta (I_{\Delta (C(v-1,2))})$. 
  
  Let $x_0$ be a new vertex and write
  $\{x_0\} \ast \Gamma(v-1)$ for the cone of $\Gamma(v-1)$
  over $x_0$.  In other words,
\[
   \{x_0\} \ast \Gamma(v-1) =
   \{ \{ x_0 \} \cup F \, : \, F \in \Gamma(v-1) \}
   \cup \Gamma(v-1).
\]
 By using the formula (\ref{cyclic}) it follows easily that
  \begin{displaymath}
    \beta_i (I_{\Delta(C(v,2))}) = \left\{
    \begin{alignedat}{3}
      &f_0 (\{x_0\} \ast \Gamma(v-1)) + v-3, &\quad &i=0, \\
      &f_i (\{x_0\} \ast \Gamma(v-1)) + \binom{v-2}{i+1}, 
      &\quad &1 \leq i \leq v-5, \\
      &f_{v-4} (\{x_0\} \ast \Gamma(v-1)) + v-3, &\quad &i=v-4, \\
      &1, &\quad &i=v-3. 
    \end{alignedat}
    \right. 
  \end{displaymath}
Let $x_1, \ldots, x_{v-3}$ be new vertices and $\Gamma'$
the simplicial complex consisting of all subsets of
$\{ x_0, x_1, \ldots, x_{v-3} \}$.  We then introduce the
simplicial complex $\Gamma(v)$ by setting
  \begin{displaymath}
    \Gamma (v) = (\{ x_0 \} \ast \Gamma (v-1)) 
      \cup \left(
      \Gamma'
                  \setminus \left\{ \{ x_0, x_1, \ldots, x_{v-3} \}, 
                                    \{ x_1, \ldots, x_{v-3} \} \right\} 
           \right). 
  \end{displaymath}
  Since $f_{v-3} (\{ x_0 \} \ast \Gamma (v-1))
  = f_{v-4} (\Gamma (v-1)) = 1$, one has 
  $\beta_i (I_{\Delta(C(v,2))}) = f_i (\Gamma (v))$ for all $i$, 
  as desired. 

  \par
  Next, let $d'>1$. 
  Again, we show that, by using induction on $v$, the Betti 
  sequence $\beta(I_{\Delta(C(v,d))})$ is the $f$-vector 
  of a simplicial complex.
  When $v=d+2$, the Betti sequence of $I_{\Delta (C(v,d))}$ 
  is $(2,1)$, which is the $f$-vector of a $1$-simplex. 
  
  Let $v > d+2$ and suppose that
  there exists a simplicial complex $\Gamma^\sharp = \Gamma (v-1,d)$ 
  such that $f (\Gamma^\sharp) = \beta (I_{\Delta(C(v-1,d))})$. 
  On the other hand, since we are working on induction on $d'$, 
  it follows that
  there exists a simplicial complex $\Gamma^\flat = \Gamma (v-2,d-2)$ 
  such that $f (\Gamma^\flat) = \beta (I_{\Delta (C(v-2,d-2))})$. 
  We will assume that the vertex set of 
  $\Gamma^\sharp$ and that of $\Gamma^\flat$ are disjoint. 

  \par
  Let $x_0$ be a new vertex.  
  Again, by using the formula (\ref{cyclic}) 
  it follows easily that
  \begin{displaymath}
    \beta_i (I_{\Delta(C(v,d))}) = \left\{
    \begin{alignedat}{3}
      &f_0 (\{ x_0 \} \ast \Gamma^\sharp) + f_0 (\Gamma^\flat) - 1, 
        &\quad &i=0, \\
      &f_i (\{ x_0 \} \ast \Gamma^\sharp) + f_i (\Gamma^\flat), 
      &\quad &1 \leq i \leq v-d-3, \\
      &f_{v-d-2} (\{ x_0 \} \ast \Gamma^\sharp) 
      + f_{v-d-2} (\Gamma^\flat) - 1, 
        &\quad &i=v-d-2, \\
      &1, &\quad &i=v-d-1. 
    \end{alignedat}
    \right. 
  \end{displaymath}
  In other words, 
  \begin{displaymath}
    \beta_i (I_{\Delta(C(v,d))}) = 
      f_i (\{ x_0 \} \ast \Gamma^\sharp) + f_i (\Gamma^\flat) - 1, 
        \quad i=0, v-d-2, v-d-1. 
  \end{displaymath}
Let $y_0$ be a vertex of $\Gamma^\flat$. 
  Let $F \in \Gamma^\flat$ be the unique face of 
  dimension $v-d-1$ and 
  $G$ a maximal proper subset of $F$. Then
  the simplicial complex 
  \begin{displaymath}
    \Gamma (v,d) = (\{ y_0 \} \ast \Gamma^\sharp) 
      \cup ( \Gamma^\flat \setminus \{ F,G \}) 
  \end{displaymath}
  satisfies 
  $\beta_i (I_{\Delta(C(v,d))}) = f_i (\Gamma (v,d))$ for all $i$,
  as desired.
  \, \, \, \, \, \, \, \, \, \, 
  \, \, \, \, \, \, \,
\end{proof}

\begin{Conjecture}
{\em
The Betti sequence arising from a pure resolution 
of a monomial ideal is
the $f$-vector of a simplicial complex.
}
\end{Conjecture}

\section{Edge ideals of chordal graphs}
Let $V$ be the vertex set and $G$ a finite graph on $V$
having no loop and no multiple edge.
Let $E(G)$ denote the edge set of $G$. 
We write $S = K[\{ x : x \in V\}]$ for the polynomial ring
in $|V|$ variables over a field $K$ 
with each $\deg x = 1$.  The {\em edge ideal} 
of $G$ is the ideal $I(G)$ of $S$ generated by
those monomials $x y$ with $\{x, y\} \in E(G)$.  

We cannot escape from the temptation to ask if
the Betti sequence of the edge ideal of a finite graph
can be the $f$-vector of a simplicial complex.  
Unfortunately, as was stated explicitly in Introduction, 
the Betti sequence of the edge ideal of the cycle 
of length $6$ cannot be the $f$-vector of
a simplicial complex.
However, it turn out to be true that
the Betti sequence of the edge ideal of a finite 
chordal graph can be the $f$-vector of a simplicial complex 
(Theorem \ref{claim:chordal_graph}).
Recall that a finite graph $G$ is {\em chordal}
if each cycle of $G$ of length $> 3$ has a chord.

\begin{Theorem}
  \label{claim:chordal_graph}
  Given an arbitrary chordal graph $G$,
  there exists a simplicial complex $\Delta$ 
  whose $f$-vector $f(\Delta)$ 
  coincides with the Betti sequence $\beta(I(G))$
of the edge ideal $I(G)$.
\end{Theorem}

\par
The recursive-type formula 
(\cite[Theorem 5.8]{Ha-Tuyl}) 
due to H\`{a} and Van Tuyl 
will be indispensable to achieve the proof of 
Theorem \ref{claim:chordal_graph}.

Let, as before, $G$ be a finite graph on  
$V$ and $E(G)$ its edge set. 
Given a subset $W \subset V$,
the \textit{restriction} $G$ to $W$
is the finite graph $G_W$ on $W$ whose edges
are those edges 
$e \in E(G)$ with $e \subset W$. 
The \textit{neighborhood} of a vertex $v$ of $G$
is the subset $N(v) \subset V$ consisting of those vertices
  $u$ of $G$ with $\{ u, v \} \in E(G)$.
We write $G \setminus e$, where $e \in E(G)$,
for the subgraph of $G$ 
which is obtained by removing $e$ from $G$. 
The \textit{distance} $\dist_G (e, e')$
of two edges $e, e' \in E(G)$  
is the smallest integer 
$\ell \geq 0$ for which there is
a sequence 
$e = e_0, e_1, \ldots, e_\ell = e'$,
where each $e_i \in E(G)$, 
with $e_{i} \cap e_{i+1} \neq \emptyset$
for all $i$.  

A {\em complete graph} on $V$ is the finite graph on $V$
such that $\{x, y\}$ is its edge 
for all $x, y \in V$ with $x \neq y$.

\begin{Lemma}[H\`{a} and Van Tuyl]
  \label{claim:rec_formula}
  Let $G$ be a chordal graph and $E(G)$ its edge set. 
  Suppose that $e = \{ u, v \}$ is an edge of $G$
  such that $G_{N(v)}$ is a complete graph.
  Let $t = |N(u) \setminus \{ v \}|$ and $G'$ the subgraph of $G$ with 
\[
    E(G') = \{ e' \in E(G) \, : \, \dist_G (e, e') \geq 3 \}. 
\]
  Then each of $G \setminus e$ and $G'$ is chordal and 
  \begin{eqnarray}
  \label{formula}
    \beta_i (I(G)) = \beta_i (I(G \setminus e)) 
                   + \sum_{\ell=0}^i \binom{t}{\ell} \beta_{i-\ell-1} (I(G'))
  \end{eqnarray}
  for all $i \geq 0$, where $\beta_{-1} (I(G')) = 1$. 
\end{Lemma}

\begin{Remark}
\label{remark}
{\em
(a) 
In Dirac \cite{Dirac}
it is proved that 
a finite graph $G$ is chordal
if and only if $G$ possesses a ``perfect elimination ordering.''
This fact guarantees the existence of a vertex $v$ of 
a chordal graph $G$
such that $N(v)$ is a complete graph.

(b) Let $N(u) = \{ v, x_1, \ldots, x_t \}$.
  Since $N(v)$ is complete,
  if $\{ v, z \} \in E(G)$, then 
  $\{ u, z \} \in E(G)$.  
  In particular,
  if $z \not\in \{ u, v, x_1, \ldots, x_t \}$, 
  then  
  $\{ v, z \} \not\in E(G)$.
  Thus
  an edge $e'$ of $G$ satisfies
  $\dist_G (e, e') \leq 2$ 
  if and only if 
  $e' \cap \{ u, v, x_1, \ldots, x_t \} \neq \emptyset$. 
  Let $W$ denote the subset of $V$ consisting of those vertices
  $z$ such that there is $e' \in E(G')$
  with $z \in e'$. 
  In particular $W \subset V \setminus \{u, v, x_1, \ldots, x_t\}$.
  Obviously $G' \subset G_W$. 
  Since none of the vertices
  $u, v, x_1, x_2, \ldots, x_t$ belongs to $W$,
  one has 
  $\dist_G (e, e') \geq 3$ for $e' \in E(G_W)$. 
  Hence $G_W \subset G'$.  Thus $G' = G_W$.   
}
\end{Remark}

\begin{Example}
{\em 
  Let $G$ be the chordal graph on $\{ y_1, \ldots, y_8 \}$
  drawn below. 
  Let $v = y_1$, $u = y_2$, and $e = \{ u, v \}$.
  Then $G_{N(v)}$ is a complete graph,
  $N (u) \setminus \{ v \} = \{ y_3, y_4, y_5 \}$,
  $t=3$ and $G' = G_{\{ y_6, y_7, y_8 \}}$.

  
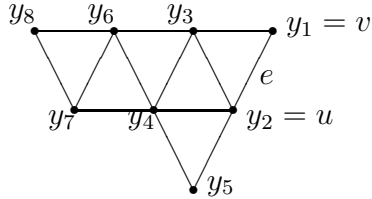
\par
  \begin{center}
    \begin{picture}(130,100)
      \thinlines
      \put(20,80){\line(1,0){90}}
      \put(35,50){\line(1,0){60}}
      \put(20,80){\line(1,-2){15}}
      \put(50,80){\line(1,-2){30}}
      \put(80,80){\line(1,-2){15}}
      \put(110,80){\line(-1,-2){30}}
      \put(80,80){\line(-1,-2){15}}
      \put(50,80){\line(-1,-2){15}}
      \put(20,80){\circle*{3}}
      \put(50,80){\circle*{3}}
      \put(80,80){\circle*{3}}
      \put(110,80){\circle*{3}}
      \put(35,50){\circle*{3}}
      \put(65,50){\circle*{3}}
      \put(95,50){\circle*{3}}
      \put(80,20){\circle*{3}}
      \put(105,60){$e$}
      \put(10,85){$y_8$}
      \put(40,85){$y_6$}
      \put(70,85){$y_3$}
      \put(115,80){$y_1 = v$}
      \put(25,45){$y_7$}
      \put(55,45){$y_4$}
      \put(100,45){$y_2 = u$}
      \put(85,20){$y_5$}
    \end{picture}
  \end{center}
%
The Betti sequences of $I(G)$, $I(G \setminus e)$ and $I(G')$
are 
  \begin{displaymath}
    \begin{aligned}
      \beta (I(G)) &= (13, 36, 47, 34, 13, 2), \\
      \beta (I(G \setminus e)) &= (12, 30, 33, 18, 4), \, \, \, \, \, 
      \beta (I(G')) = (3, 2). 
    \end{aligned}
  \end{displaymath}
We can easily check that these Betti sequences 
  satisfy the formula (\ref{formula}) due to H\`{a} and Van Tuyl.  
  For example, since   
  $47 = 33 + 2 \cdot \binom{3}{0} + 3 \cdot \binom{3}{1} 
  + 1 \cdot \binom{3}{2}$,
  one has
  \begin{displaymath}
    \begin{aligned}
     \beta_2 (I(G)) =
     \beta_2 (I(G \setminus e)) 
       + \binom{3}{0} \beta_{1} (I(G')) + \binom{3}{1} \beta_{0} (I(G')) 
       + \binom{3}{2} \beta_{-1} (I(G')).
    \end{aligned}
  \end{displaymath}
}
\end{Example}

\begin{Lemma}
\label{Hochster}
Let $G$ be an arbitrary graph on $V = V(G)$
and $W$ a subset of $V$.  Then one has
\[
\beta_i (I(G)) \geq \beta_i (I(G_W))
\]
for all $i$.
\end{Lemma}

\begin{proof}
  Since $I(G)$ and $I(G_W)$ are squarefree monomial ideals, 
  there exist simplicial complexes $\Delta$ on $V$
  and ${\Delta}'$ on $W$
  such that $I_{\Delta} = I(G)$ 
  and $I_{{\Delta}'} = I(G_W)$. 
  Hochster's formula  
  \cite[Corollary 4.9, p. 64]{Stanley} says that 
  \begin{displaymath}
    \begin{aligned}
      {\beta}_i (I(G)) &= {\beta}_i (I_{\Delta}) 
      = \sum_{U \subset V} 
          \dim_K \widetilde{H}_{|U| - i - 2} ({\Delta}_U; K), \\
      {\beta}_i (I(G_W)) &= {\beta}_i (I_{{\Delta}'}) 
      = \sum_{U' \subset W} 
          \dim_K \widetilde{H}_{|U'| - i - 2} ({\Delta}_U'; K). 
    \end{aligned}
  \end{displaymath}
  What we must prove is that ${\Delta}_U = {\Delta}'_U$ 
  whenever $U \subset W$. 

  \par
  Let $F \in {\Delta}_U$.
  Then, for all $\{ x, y \} \subset F$,
  one has $\{ x, y \} \not\in E(G)$.
  In particular $\{ x, y \} \not\in E(G_W)$. 
  Thus $F \in {\Delta}'$ and $F \in {\Delta}'_U$. 
  Conversely, let $F \in {\Delta}'_U$. 
  Then, for all $\{ x, y \} \subset F$, 
  one has $\{ x, y \} \not\in E(G_W)$. 
  Since $\{ x, y \} \subset F \subset U \subset W$, one has 
  $\{ x, y \} \not\in E(G)$. 
  Hence $F \in {\Delta}$ and $F \in {\Delta}_U$, as desired. 
  \, \, \, \, \, \, \, \, \, \, 
  \, \, \, \, \, \, 
\end{proof}

\begin{Lemma}
  \label{claim:colon_ideal}
  Let $S$ be a polynomial ring over a field $K$. 
  
 {\em (a)}
  Let $I \subset S$ be a squarefree monomial ideal and 
    $x$ a variable of $S$. Then 
    \begin{displaymath}
      {\beta}_i (I) \geq {\beta}_i (I:x)
    \end{displaymath}
    for all $i$.
    
  {\em (b)}
  Let $I$ and $J$ be monomial ideals of $S$
  and $G(I)$ (resp. $G(J)$) the minimal system of
  monomial generators of $I$ (resp. $J$).
  Suppose that $\supp(u) \cap \supp(v) = \emptyset$
  for all $u \in G(I)$ and for all $v \in G(J)$, where
  $\supp(u)$ is the set of variables $x \in V$
  which divides $u$.
    Then, for all $i$, one has 
    \begin{displaymath}
      \beta_{i} (S/(I+J)) 
        = \sum_{m=0}^{i} \beta_{i-m} (S/I) \beta_{m} (S/J). 
    \end{displaymath}
\end{Lemma}

\begin{proof}
  (a) Let $S = K [x_1, x_2, \ldots, x_n]$, $x=x_1$, 
  and $R = K[x_2, \ldots, x_n]$. 
  Since the variable $x$ does not appear in the minimal system of 
  monomial generators of $I:x$,
  it follows that
  $J = (I:x) \cap R$ has the same minimal set of 
  monomial generators as that of $I:x$. 
  Hence $\beta^R_i (J) = \beta^S_i (I:x)$,
  where $\beta^R_i (J)$ is the $i$th total Betti number
  of $J \subset R$ and $\beta^S_i (I:x)$ 
  is that of $I:x \subset S$.
  We claim $\beta^S_i (I) \geq \beta^R_i(J)$. 
  
  \par
  Let $F_{\bullet}$ be a minimal graded free resolution of $S/I$ on $S$.
  Since $x - 1$ 
  is a non-zero divisor of $S/I$, by using 
  \cite[Proposition 1.1.5]{BH},
  it follows that 
  $F_{\bullet} \otimes S/(x-1)$ is a free resolution of 
  $S/I \otimes S/(x-1)$ on $S \otimes S/(x-1)$.
  Since 
  \begin{displaymath}
    R/J \cong S/I \otimes_S S/(x-1), 
    \qquad R \cong S/(x-1) \cong S \otimes S/(x-1), 
  \end{displaymath}
  $F_{\bullet} \otimes S/(x-1)$ is a free resolution of 
  $R/J$ on $R$.  

\smallskip

  \par
  (b) Let $F_{\bullet}$ (resp.\  $G_{\bullet}$) be a minimal graded free 
  resolution of $S/I$ (resp.\  $S/J$) 
  and $T_{\bullet}^{(I)}$ (resp.\  $T_{\bullet}^{(J)}$) the Taylor resolution 
  of $S/I$ (resp.\  $S/J$). 
  Then $T_{\bullet}^{(I)} \otimes T_{\bullet}^{(J)}$ is isomorphic to 
  the Taylor resolution of $S/(I+J)$. Thus 
  \begin{displaymath}
    H_i (F_{\bullet} \otimes G_{\bullet}) 
    \cong \Tor_i^S (S/I, S/J) 
    \cong H_i (T_{\bullet}^{(I)} \otimes T_{\bullet}^{(J)}) = 0. 
  \end{displaymath}
  Hence $F_{\bullet} \otimes G_{\bullet}$ is a graded free resolution of 
  $S/I \otimes S/J \cong S/(I+J)$. 
  In particular
  $F_{\bullet} \otimes G_{\bullet}$ is minimal. 
 \, \, \, \, \, \, \, \, \, \,  
 \, \, \, \, \, \, \, \, \, \, 
 \, \, \, \, \, \, \, \, \, \, 
 \, \, \, \, \, \, \, \, \, \, 
 \, 
\end{proof}

 
\begin{Lemma}
  \label{claim:betti_cone}
  Let $G$ be an arbitrary graph on $V$ 
  and let $W$ be a subset of $V$. 
  Suppose that $G_{V \setminus W}$ contains edges 
  \begin{displaymath}
    \{ u, x_1 \}, \{ u, x_2 \}, \ldots, \{ u, x_t \}, 
  \end{displaymath}
  where $t \geq 1$ is an integer and 
  where $u, x_1, x_2, \ldots, x_t$ are distinct vertices of $G$. 
  If $\{ u, z \} \not\in E(G)$ for all $z \in W$, then 
  \begin{displaymath}
    \beta_i (I(G)) \geq \sum_{m=0}^{i+1} \binom{t}{m} \beta_{i-m} (I(G_W)), 
  \end{displaymath}
  for all $i \geq 0$, where $\beta_{-1} (I(G_W)) = -1$. 
\end{Lemma}

\begin{proof}
  Set $V' = \{ u, x_1, \ldots, x_t \}$. 
  Lemma \ref{Hochster}
  together with
  Lemma \ref{claim:colon_ideal} (a) says that 
  \begin{displaymath}
    \beta_i (I(G)) \geq \beta_i (I(G_{V' \cup W})) 
    \geq \beta_i (I(G_{V' \cup W}) : u). 
  \end{displaymath}
  Since $\{ u, z \} \not\in E(G)$ 
  for all $z \in W$, it follows that
  \begin{displaymath}
    I(G_{V' \cup W}) : u = I(G_W) + (x_1, x_2, \ldots, x_t). 
  \end{displaymath}
  Then, since $V' \cap W = \emptyset$,
  by using Lemma \ref{claim:colon_ideal} (b), 
  one has 
  \begin{displaymath}
    \begin{aligned}
      \beta_i (I(G_W) + (x_1, x_2, \ldots, x_t)) 
      &= \beta_{i+1} (S/(I(G_W) + (x_1, x_2, \ldots, x_t))) \\
      &= \sum_{m=0}^{i+1} 
           \beta_{i+1-m} (S/I(G_W)) \beta_{m} (S/(x_1, x_2, \ldots, x_t)) \\
      &= \sum_{m=0}^{i+1} 
           \binom{t}{m} \beta_{i-m} (I(G_W)), 
    \end{aligned}
  \end{displaymath}
  as required.
\, \, \, \, \, \, \, \, \, \, 
\, \, \, \, \, \, \, \, \, \, 
\, \, \, \, \, \, \, \, \, \, 
\, \, \, \, \, \, \, \, \, \,  
\, \, \, \, \, \, \, 
  \end{proof}

Let $\Delta$ be a simplicial complex on the vertex set
$V$ and let $x$ be a new vertex. 
The {\em cone} of $\Delta$ over $x$ is 
the simplicial complex 
\begin{displaymath}
  \cone (\Delta) = \{ \{ x \} \cup F \, : \, F \in \Delta \} \cup \Delta 
\end{displaymath}
on $V \cup \{ x \}$.
Moreover, by setting $\cone^{0} (\Delta) = \Delta$,
the $t$th cone of $\Delta$ 
is defined recursively by 
\begin{displaymath}
  \cone^t ({\Delta}) = \cone (\cone^{t-1} (\Delta)).
\end{displaymath}
It follows that
\[
f_i(\cone^t ({\Delta})) = 
\sum_{\ell=0}^{i+1} \binom{t}{\ell} f_{i-\ell} ({\Delta})
\]
for all $i$.

We are now in the position to give a proof of 
Theorem \ref{claim:chordal_graph}.
Recall that the Stanley--Reisner ideal
$I_\Delta \subset S$ is {\em squarefree lexsegment} 
(\cite{AHH}) 
if, 
for all monomials
$u$ and $v$ of $S$ with $\deg u = \deg v$ 
and with $v <_{\lex} u$
such that $v \in I_{\Delta}$, 
one has $u \in I_{\Delta}$,
where $<_{\lex}$ is the lexicographic order
induced by a (fixed) ordering of the variables 
of $S$. 
Given a simplicial complex $\Delta$, 
there is a unique simplicial complex $\Delta^{\lex}$
such that $I_{\Delta^{\lex}}$ is 
squarefree lexsegment
with $f(\Delta) = f(\Delta^{\lex})$.

\begin{proof}[Proof of Theorem \ref{claim:chordal_graph}]
  Our proof will proceed  
  by using induction on the number of edges of $G$. 
  If $G$ possesses only one edge $\{ x, y \}$, then
  $I(G) = (xy)$ and 
  \begin{displaymath}
    \beta_i (I(G)) = \left\{
    \begin{aligned}[3]
      &1, &\quad &i = 0, \\
      &0, &\quad &i \neq 0. 
    \end{aligned}
    \right. 
  \end{displaymath}
  Thus its Betti sequence is equal to 
  the $f$-vector of a $0$-simplex. 

  \par
  Now, suppose that $G$ possesses at least two edges
  and that, for an arbitrary chordal graph $\Gamma$ 
  with $|E(\Gamma)| < |E(G)|$, the Betti sequence
  $\beta(I (\Gamma))$ is the $f$-vector $f(\Delta_\Gamma)$
  of a simplicial complex
  $\Delta_\Gamma$.

  \par
  Let $e = \{u, v \}$ be an edge of $G$ such that 
  $G_{N(v)}$ is complete. 
  Work with the same notation as in Lemma \ref{claim:rec_formula}
  and in Remark \ref{remark} (b). 
  One has 
  \begin{equation*}
    \label{eq:rec_formula}
    \beta_i (I(G)) = \beta_i (I(G \setminus e)) 
        + \sum_{\ell=0}^i \binom{t}{\ell} \beta_{i-\ell-1} (I(G_W)),
  \end{equation*}
  Since each of $G \setminus e$ and $G_W$ is a subgraph of $G$
  with $e \not\in E(G \setminus e)$ and $e \not\in E(G_W)$, 
  the hypothesis of induction guarantees the existence 
  of simplicial complexes 
  $\Delta_{G \setminus e}$ and $\Delta_{G_W}$ such that 
  \begin{equation*}
    \label{eq:f=b}
    f_i (\Delta_{G \setminus e}) = \beta_i (I(G \setminus e)), 
    \quad
    f_i (\Delta_{G_W}) = \beta_i (I(G_W)). 
  \end{equation*}
Thus what we must prove is the existence of a simplicial complex
$\Delta$ with
  \begin{equation}
    \label{final}
f_i (\Delta) = f_i ({\Delta}_{G \setminus E}) 
    + \sum_{\ell=0}^{i} \binom{t}{\ell} f_{i-\ell-1} ({\Delta}_{G_W}).
  \end{equation}

It follows from Lemma \ref{claim:betti_cone} that
  \begin{equation*}
    \label{eq:WTSineq}
    \beta_i (I (G \setminus e)) 
    \geq \sum_{m=0}^{i+1} \binom{t}{m} \beta_{i-m} (I(G_W)). 
  \end{equation*}
In other words, 
  \begin{displaymath}
    f_i (\Delta_{G \setminus e}) 
    \geq \sum_{m=0}^{i+1} \binom{t}{m} f_{i-m} (\Delta_{G_W})
    = f_i (\cone^t ({\Delta}_{G_W})). 
  \end{displaymath}
Thus, by choosing $\Delta_{G \setminus e}$
for which $I_{\Delta_{G \setminus e}}$ 
is squarefree lexsegment,
we assume that $\Delta_{G \setminus e}$
contains a subcomplex $\Delta'$ whose
$f$-vector coincides with that of 
$\cone^t ({\Delta}_{G_W})$.

We introduce the simplicial complex $\Delta$ by setting 
 \begin{displaymath}
    \Delta = {\Delta}_{G \setminus E} \cup \cone ({\Delta}'), 
  \end{displaymath}
  where the new vertex of $\cone ({\Delta}')$ cannot be
  a vertex of ${\Delta}_{G \setminus E}$.
Then 
  \begin{displaymath}
    \begin{aligned}
    f_i (\Delta) - f_i ({\Delta}_{G \setminus E}) 
    &= f_i (\cone ({\Delta}')) - f_i ({\Delta}') \\
    &= f_{i-1}({\Delta}') \\
    &= f_{i-1}({\cone}^t (\Delta_{G_W})) \\
    &= \sum_{\ell=0}^{i} \binom{t}{\ell} f_{i-\ell-1} ({\Delta}_{G_W}). 
    \end{aligned}
  \end{displaymath}
Thus the simplicial complex satisfies the equality 
(\ref{final}), as desired. 
\, \, \, \, \, \, \, \, \, 
\end{proof}

\section{Gorenstein Monomial Ideals}
We now turn to the discussion on Betti sequences of 
Gorenstein monomial ideals.
Let, as before, $S = K[x_1,\dots,x_n]$
denote the polynomial ring in $n$ variables 
over a field $K$ with each $\deg x_i = 1$.
Recall that 
a homogeneous ideal $I \subset S$ is Gorenstein if $S/I$
is a Gorenstein ring.
If $I \subset S$ is Gorenstein,
then its Betti sequence 
$\beta(I)=(\beta_0(I), \beta_1(I), \ldots, \beta_{p}(I))$ 
is symmetric, that is, 
$\beta_i(I)=\beta_{p-1-i}(I)$ for all $i$,
where $p=\projdim(I)$ and where $\beta_{-1}(I) = 1$. 

Let $I \subset S$ be a Gorenstein monomial ideal with 
$\projdim(I)=p$.
If $p=1$, then $\beta(I)=(2,1)$ by the Hilbert--Burch theorem 
\cite[Theorem 1.4.17]{BH}.
If $p=2$, then there exists an odd integer $m \geq 3$ such that
$\beta(I)=(m,m,1)$ by the structure theorem 
due to Buchsbaum and Eisenbud (\cite[Theorem 3.4.1]{BH}).
In fact, 
these facts characterize the Betti numbers of 
Gorenstein (monomial) ideals with $\projdim(I) \leq 2$.
For example, 
to prove the sufficiency, 
let $I$ be the Stanley--Reisner ideal
of the boundary complex of the cyclic $2m$-polytope 
with $2m+3$ vertices.  Then $I$ is a Gorenstein ideal with 
$\beta(I)=(2m+3,2m+3,1)$ for all $m \geq 1$ 
by the formula (\ref{cyclic}).

Let $p=3$.
Let $I \subset S$ be a Gorenstein monomial ideal with 
$\projdim(I) = 3$.
Since  
$(\beta_{-1}(I), \beta_{0}(I), 
\beta_{1}(I), \beta_{2}(I), \beta_{3}(I))$,
where $\beta_{-1}(I)=1$, 
is symmetric 
and since 
$\sum_{i=-1}^{3} ( - 1 )^i \beta_i(I) = 0$,
it follows that 
there exists an integer $m$ such that
$\beta(I)=(m+1,2m,m+1,1)$.
Since $I$ is a monomial ideal, 
the Taylor resolution of $I$ says that
$m= \beta_0(I) -1 \geq \projdim(I)=3$ .
Since $(m,m,1)$ is the $f$-vector of a simplicial complex 
for $m \geq 3$,
it follows from Lemma \ref{Kalaiacyclic} that
$\beta(I)$ is the $f$-vector of an acyclic simplicial complex.

\begin{Example}
\label{ex1}
{\em
Let $I=(x_1x_4,x_1x_5,x_2x_6,x_3x_7,x_4x_6,x_4x_7,x_2x_3x_5)$.
Then $I$ is Gorenstein and 
$\beta(I)=(7,12,7,1)=(6+1,2 \times 6, 6+1,1)$.
}
\end{Example}

More precisely, we can characterize the Betti numbers
of Gorenstein monomial ideals $I$ with 
$\projdim(I) = 3$.
Recall that
a monomial ideal $I \subset S$ is \textit{strongly stable}
if, for all monomials $u  \in I$ and for all $j < i$ such that $x_i$
divides $u$, 
one has $u x_j/ x_i \in I$.

\begin{Theorem}
\label{gormonomial}
Let $\beta=(m+1,2m,m+1,1)$, where $m$ is an integer with $m \geq 3$.
Then there exists a Gorenstein monomial ideal $I$ of a polynomial
ring with
$\beta(I)=\beta$ if and only if $m=3$ or $m \geq 5$.
\end{Theorem}

\begin{proof}
{\bf (``If'')}
Let $m \geq 3$ be odd.
Then there exists a Gorenstein monomial ideal $J \subset S$ 
with
$\beta(J)=(m,m,1)$. Let $y$ be a new variable and $S' = S[y]$.
Then the ideal $I=J+(y)$ is a Gorenstein monomial ideal 
with $\projdim(I) = 3$ and $\beta_0(I)=m+1$.

Let $m \geq 6$ be even.
Example \ref{ex1} yields an example of $m = 6$.
Now, let $m = 2k + 6 \geq 8$ be even.
Given a strongly stable ideal $J \subset R=K[x_1,\dots,x_p]$ 
such that
$R/J$ is of finite length,
it follows from \cite[Theorem 9.6]{MN} and \cite[Theorem 5.3]{Mu}
that there exists a Gorenstein squarefree monomial ideal $I_{(J)}$ 
for which
$\beta_i(S/I_{(J)})=\beta_i(R/J)+ \beta_{p+1-i}(R/J)$ for all $i$.
Let $J$ be the strongly stable ideal
$$J=(x_1^2,x_1x_2,x_1x_3,x_2^{k+1},x_2^kx_3,\cdots,x_3^{k+1}) \subset K[x_1,x_2,x_3].$$
Eliahou--Kervaire formula
says that 
$\beta_0(I_{(J)})=\beta_0(J)+\beta_2(J)=2k+7$, as required.

\smallskip

{\bf (``Only If'')} 
We show, in general, that 
if $I \subset S$ is a Gorenstein monomial ideal with 
$\projdim(I)=p-1 \geq 3$, then  
$\beta_0(I) \ne p+1$.
Suppose, on the contrary,
that there exists a Gorenstein monomial ideal with
$\projdim(I)=p-1 \geq 3$ and $\beta_0(I)= p+1$.
By taking the polarization (\cite[Lemma 4.2.16]{BH}) of $I$,
we assume that $I$ is squarefree.
Without loss of generality,
we assume that no variable $x_i$ is contained in $I$.
Let $\Delta$ (resp. $\Delta'$)
denote the simplicial complex whose Stanley--Reisner ideal is 
$I$ (resp. the  $I:x_i$).  
Then $\Delta'$ is the star (\cite[Definition 5.3.4]{BH}) 
of $\Delta$ of the face $\{ i \}$.
Hence $I:x_i$ is a Gorenstein ideal 
with $\dim(S/I)=\dim(S/(I:x_i))$. 
In particular $\projdim(I)=\projdim(I:x_i)$. 
Thus, in case of $\beta_0(I) = \beta_0(I:x_i)$, 
we replace $I$ with $I:x_i$.
Hence, for each variable $x_k$ 
which appears in the minimal system of monomial generators 
of $I$, we assume that $\beta_0(I:x_k) <p+1$.
On the other hand, since 
$\projdim(I)=\projdim(I:x_k) \leq \beta_0(I:x_k) - 1$,
it follows that $\beta_0(I:x_k) = p$ and
$I:x_k$ is a complete intersection.

Let $G(I)=\{u_1,\dots,u_{p+1}\}$ 
be the minimal system of monomial generators of $I$.
Say, $x_1$ divides $u_1$
and, since $\beta_0(I:x_1)=p$,  
$u_1/x_1$ divides $u_{p+1}$.
Let $u_1=x_1 x_F$ and $u_{p+1}=x_F x_G$,
where $x_F = \prod_{i \in F} x_i$ with 
$F \subset [n]$.
Then
\[
I:x_1=(\tilde u_1,\tilde u_2,\dots,\tilde u_p),
\]
where $\tilde u_k=u_k/x_1$ (resp. $\tilde u_k = u_k$)
if $x_1$ divides (resp. does not divide) $u_k$. 
In particular $\tilde u_1 = x_F$.
Since $I:x_1$ is a complete intersection,
it follows that
\begin{eqnarray}
\label{siki1}
\supp(\tilde u_s) \cap \supp(\tilde u_t) = \emptyset
\end{eqnarray}
if $s \neq t$, where $\supp(\tilde u_s)$ stands for
the set of variables $x_k$ which divides $\tilde u_s$. 
If there is $2 \leq k \leq p$ 
with $u_k = \tilde u_k$,
then,
since $\beta_0(I:x_j)=p$
for all $x_j \in \supp(u_k)$,
it follows from (\ref{siki1}) that
$u_k$ must divide $u_{p+1}$, a contradiction.
Thus $\tilde u_k = u_k / x_1$ for each $1 \leq k \leq p$.
Let $j \in F$.
Then, by (\ref{siki1}), $x_j \not\in \supp(u_k)$
for $k = 2, \ldots, p$.
Since $\beta_0(I:x_j)=p$,
there is $k$ with $2 \leq k \leq p$
such that
either $u_1/x_j$ or $u_{p+1}/x_j$ must divide $u_k$.
If $u_1/x_j$ divides $u_k$, then 
$u_1 = x_1x_j$ by (\ref{siki1}).
Thus 
$\beta_0(I:x_j)=p=2$,
a contradiction.
If $u_{p+1}/x_j=x_Gx_F/x_j$ divides $u_k$,
then, again by (\ref{siki1}),
one has $u_1 = x_1 x_F = x_1 x_j$ 
and $p=2$, a contradiction.
\, \, \, \, \, \, \, \, \, \, 
\, \, \, \, \, \, \, \, \, \, 
\, \, \, 
\end{proof}

The technique appearing in the ``If'' part of the proof 
of Theorem \ref{gormonomial} together with the result shown
in the ``Only If'' part of Theorem \ref{gormonomial}
yields the following
 
\begin{Corollary}
Fix integers $m \geq 4$ and $p \geq 3$.
Then there exists a Gorenstein monomial ideal $I$ of
$K[x_1, \ldots, x_n]$, where $n$ is enough large,
with $\beta_0(I) = m$ 
and $\projdim(I) = p$ if and only if  
$m \geq p + 1$ with $m \ne p+2$.
\end{Corollary}

\end{document}